\documentclass[12pt]{amsart}
\usepackage{amsmath,amssymb, xypic,verbatim,amscd,color}
\usepackage{graphicx}
\usepackage{colordvi}
\usepackage{times}

\numberwithin{equation}{section}

\newcommand\Pdescend{P}

\newcommand\Spin{\operatorname{Spin}}
\newcommand\SO{\operatorname{SO}}

\newcommand\point{x}

\newcommand\Sp{\operatorname{Sp}}

\newcommand\Pic{\operatorname{Pic}}

\newcommand\tensor{\otimes}
\newcommand\ml{\mathcal{L}}

\newcommand\reg{\operatorname{reg}}

\newcommand\mg{\mathcal{G}}

\newcommand\mq{\mathcal{Q}}

\newcommand\SL{\operatorname{SL}}

\newcommand{\leto}[1]{\stackrel{#1}{\to}}

\newcommand\PP{\mathcal{P}}

\newcommand\pii{p}

\newtheorem{theorem}{Theorem}[section]
\newtheorem{remark}[theorem]{ Remark}

\newtheorem{corollary}[theorem]{Corollary}

\newtheorem{proposition}[theorem]{Proposition}

\newtheorem{lemma}[theorem]{Lemma}

\newtheorem{defi}[theorem]{Definition}

\begin{document}
\title[Theta characteristics and strange duality ]{Orthogonal bundles, theta characteristics and the symplectic strange duality}
\author{Prakash Belkale}\thanks{The author was partially
supported by the NSF} \maketitle

\begin{abstract} A basis for the space of generalized theta functions of level one for the spin
groups, parameterized by the theta characteristics on a curve, is
shown to be projectively flat over the moduli space of curves (for
Hitchin's connection). The symplectic strange duality conjecture,
conjectured by Beauville is shown to hold for all curves of genus
$\geq 2$, by using Abe's proof of  the conjecture for generic
curves, and the above monodromy result.
\end{abstract}

\section{Introduction}
Consider the moduli stacks $\mathcal{M}_{\Spin(r)}(X)$ and
$\mathcal{M}_{\SO(r)}(X)$  of principal  $\Spin(r)$, and
$\SO(r)$-bundles, $r\geq 3$ on a smooth connected projective curve
$X$ of genus $g\geq 2$ over $\Bbb{C}$.
 Let $\mathcal{M}_{\SO(r)}(0)$ be the connected component of $\mathcal{M}_{\SO(r)}(X)$, which contains
 the trivial $\SO(r)$-bundle.
 There is a natural map
$$\pii:\mathcal{M}_{\Spin(r)}\to\mathcal{M}_{\SO(r)}(0).$$

A line bundle $\kappa$ on $X$ is said to be a theta characteristic
if $\kappa^{\otimes 2}$ is isomorphic to the canonical bundle $K_X$.
The set of theta characteristics $\theta(X)$ forms a torsor for the
$2$-torsion $J_2(X)$ in the Jacobian of $X$, and hence $|\theta(X)|=
2^{2g}$. Recall that a theta characteristic $\kappa$ is said to be
even (resp. odd) if $h^0(\kappa)$ is even (resp. odd).

 For each theta-characteristic $\kappa$ on
$X$ there is a line bundle $\mathcal{P}_k$ on $\mathcal{M}_{\SO(r)}$
with a canonical section $s_{\kappa}$
 (see the
 pfaffian construction in \cite{ls,bls}). On
$\mathcal{M}_{\SO(r)}(0)$, $s_{\kappa}=0$ if and only if  both
$\kappa$ and $r$ are odd.

For theta characteristics $\kappa$ and $\kappa'$, the line
 bundle $\pii^{*}\mathcal{P}_{\kappa}$ is
isomorphic to $\pii^{*}\mathcal{P}_{\kappa'}$ (see ~\cite{ls}). Set
$\mathcal{P}=\pii^{*}\mathcal{P}_{\kappa}$ which is well defined
upto isomorphism. The line bundle $\mathcal{P}$ is the positive
generator of the Picard group of the stack $\mathcal{M}_{\Spin(r)}$.
It is known that  $\mathcal{P}$ does not descend to the moduli-space
$M_{\Spin(r)}$, (similarly $\PP_{\kappa}$ does not descend to the
moduli-space $M_{\SO(r)}$). Clearly, $\PP$ comes equipped with
sections $s_{\kappa}$ for each theta characteristic $\kappa$, coming
from the identification $\pii^{*}\mathcal{P}_{\kappa}\leto{\sim}
\mathcal{P}$ ($s_{\kappa}$ are well defined up to scalars).

Let $\pi:\mathcal{X}\to S$ be a smooth projective relative curve of
genus $g$. Assume by passing to an \'{e}tale cover that the sheaf of
theta-characteristics on the fibers of $\pi$ is trivialized. For
$s\in S$, let $X_s=\pi^{-1}(s)$. It is known that the spaces
$H^0(\mathcal{M}_{\Spin(r)}(X_s),\mathcal{P})$ form the fibers of a
vector bundle on $S$, which is equipped with a projectively flat
connection (WZW or equivalently Hitchin's connection).

\begin{theorem}\label{phone}
For even $r$, each  section $s_{\kappa}\in
H^0(\mathcal{M}_{\Spin(r)}({X}_s),\mathcal{P})$, for
$\kappa\in\theta(X_s)$ is projectively flat.
\end{theorem}
\begin{theorem}\label{phone2}
For odd $r$, each  section $s_{\kappa}\in
H^0(\mathcal{M}_{\Spin(r)}({X}_s),\mathcal{P})$, for even
$\kappa\in\theta(X_s)$ is projectively flat.
\end{theorem}
It is known (see ~\cite{ox}) that the dimension of the space
$H^0(\mathcal{M}_{\Spin(r)}({X}_s), \mathcal{P})$ is equal to the
number of theta characteristics (if $r$ is odd, the number of even
theta characteristics). It has been proved by Pauly and Ramanan (see
Proposition 8.2 in ~\cite{PR}), that in Theorems ~\ref{phone} and
~\ref{phone2}, the sections are linearly independent, and hence form
a basis. Our methods give a new proof of this result of Pauly and
Ramanan.

We use Theorem ~\ref{phone} to show that the symplectic strange
duality formulated by Beauville ~\cite{bea} is, in a suitable sense,
projectively flat: Hence it is an isomorphism for all curves of a
given genus if it is an isomorphism for some curve of that genus
(see Corollary ~\ref{dual}). Takeshi Abe ~\cite{a1,a2} has recently
formulated a very interesting  parabolic generalization of
Beauville's conjecture, and has proved this conjecture for generic
curves by using powerful degeneration arguments. His results imply
Beauville's conjecture for generic curves\footnote{ Abe's parabolic
generalization is essential in his proof of Beauville's conjecture
for generic curves by degeneration.}. Therefore Abe's results
(together with Corollary ~\ref{dual}) imply that the symplectic
strange duality conjecture of Beauville holds for all curves. It
should be pointed out that Abe's parabolic symplectic duality
conjecture has { not} yet been shown to hold for all curves.

We would like to point out that Theorem ~\ref{phone} and
~\ref{phone2} {\it do not} imply that the global projective Hitchin
monodromy on the vector spaces $H^0(\mathcal{M}_{\Spin(r)}({X}_s),
\mathcal{P})$ is finite. The analogous question for the symplectic
group is also not known (but see Section ~\ref{etrange}).

The proofs of Theorem ~\ref{phone} and Theorem ~\ref{phone2} have
the following main ingredients.
\begin{enumerate}
\item The map $p$ can perhaps be interpreted as a ``stacky'' torsor for
$J_2(X_s)$. We will instead work over the regularly stable locus in
${M}_{\SO(r)}(0)$, over which $p$ is a torsor (using results in
~\cite{bls}).
\item By Proposition 5.2 in ~\cite{bls}, for different theta characteristics $\kappa$ and $\kappa'$, the bundles $\mathcal{P}_k$ on
 and $\PP_{\kappa'}$ on $\mathcal{M}_{\SO(r)}(X_s)$ are  not isomorphic. The isomorphism class of $\PP_{\kappa}\tensor
\PP_{\kappa'}^{-1}$ is explicitly computed in ~\cite{bls}, and this
computation  constitutes  the heart of the matter in the proofs of
Theorems ~\ref{phone} and ~\ref{phone2}.
\end{enumerate}
Avoiding technicalities, it is easy to summarize the proof of
Theorem ~\ref{phone}: Fix a theta characteristic $\kappa$. There is
an action of $J_2(X_s)$ on $(\mathcal{M}_{\Spin(r)}({X}_s),
\mathcal{P})$ which lies over a trivial  action on the pair
$(\mathcal{M}_{\SO(r)}({X}_s)(0), \mathcal{P}_{\kappa})$. Since this
action preserves the so-called geometric Segal-Sugawara tensor (see
Section ~\ref{sugawara}), it preserves Hitchin's connection on the
spaces $H^0(\mathcal{M}_{\Spin(r)}({X}_s), \mathcal{P})$. Therefore
the connection preserves  each $J_2(X_s)$-isotypical subspace of
$H^0(\mathcal{M}_{\Spin(r)}({X}_s),\mathcal{P})$. Each isotypical
subspace will be shown to contain a pfaffian section $s_{\kappa'}$.
Counting dimensions, we are then able to conclude the proof.

We will use the language of moduli-spaces and not of stacks (except
in recalling some results from ~\cite{bls}). The main technique is
to work over the regularly stable locus in ${M}_{\SO(r)}$, and to
use results of Beauville, Laszlo and Sorger ~\cite{v2,ls,bls}.
\section{Reformulation in terms of Moduli spaces}

We will use the notation, setup and results from Section 13 of
~\cite{bls}, which we recall for the benefit of the reader. Let $G$
be a simple (not necessarily simply connected) algebraic group. Let
$M_G$ denote Ramanathan's moduli space of principal semistable
$G-$bundles on a smooth projective and connected curve $X$ of genus
$g\geq 2$. Let us assume that $G$
 does not  map to $\operatorname{PGL}_2$, or that $g>2$.

\begin{defi}A  $G-$bundle on $X$ is regularly stable if it is stable
and its automorphism group is equal to the center $Z(G)$ of $G$.
\end{defi}

The open subset $M_G^{\reg}\subset M_G$ is smooth, and as pointed
out in ~\cite{bls}, the method of proof of a theorem of Faltings
(Theorem II.6 in ~\cite{fa}) implies that the  complement of
$M_G^{\reg}$ in $M_G$ is of codimension $\geq 2$.

Let $A$ be the group of principal $A'$-bundles where $A'$ is the
kernel of $\Spin(r)\to \SO(r)$ (clearly $A$ is isomorphic to $J_2$).
Denote as usual the group of one dimensional characters of $A$ by
$\hat{A}$.

Let $M_{\SO(r)}(0)$ denote the connected component of $M_{\SO(r)}$
which contains the trivial $\SO(r)$-bundle. By a result of
Beauville-Laszlo-Sorger (see the proof of Proposition 13.5 in
~\cite{bls}), the natural finite Galois covering with Galois group
$A$
$$\pii:M_{\Spin(r)}\to M_{\SO(r)}(0)$$
is  \'etale over $Y=M_{\SO(r)}^{\reg}(0)$. Set
$\widetilde{Y}=\pii^{-1}(Y)$. It follows from the proof of
Proposition 13.5 in ~\cite{bls}, that $\widetilde{Y}\subseteq
M^{\reg}_{\Spin_r}$.

Since $M_{\SO(r)}(0)-Y$ has codimension $\geq 2$ and
$\pii:M_{\Spin(r)}\to M_{\SO(r)}(0)$ is finite and dominant,
$M_{\Spin(r)}-\widetilde{Y}$ has codimension $\geq 2$. Therefore
(note that moduli spaces constructed using geometric invariant
theory are normal)
$$H^0(Y,\mathcal{O}_Y)=H^0(\tilde{Y},\mathcal{O}_{\tilde{Y}})=\Bbb{C}.$$

It is easy to see that there is a decomposition as sheaves of
$A$-modules:
$$(\pii\mid \widetilde{Y})_* \mathcal{O}_{\widetilde{Y}}=\oplus_{\chi\in \hat{A}} L_{\chi},\ L_{\chi}\in
\operatorname{Pic}(Y).$$ where as a sheaf,
$$L_{\chi}(U)=\{s\in \pii_*{\mathcal{O}}(U)\mid gs=\chi(g)s, \forall g\in A\}.$$
It is easy to verify that
\begin{itemize}
\item $H^0(Y,L_{\chi})=0$ unless $\chi=1$, $H^0(Y,L_{1})=\Bbb{C}$.
\item $p^*L_{\chi}=\mathcal{O}_{\widetilde{Y}}$,
\item $L_{\chi}\tensor L_{\chi'}=L_{\chi\chi'}$
\item  $L_{\chi}$ is not isomorphic to $L_{\chi'}$ for $\chi\neq
\chi'$.
\end{itemize}
According to Proposition 9.5 in ~\cite{ls}, the line bundle
$\PP_{\kappa}$ on $\mathcal{M}_{\SO(r)}(0)$ descends to
$M_{\SO(r)}^{\reg}(0)$, to a line bundle which we will denote by
$P_{\kappa}$ similarly the line bundle $p^*\PP_{\kappa}$ on the
moduli stack $\mathcal{M}_{\Spin(r)}$ descends to the moduli space
${M}_{\Spin(r)}^{\reg}$.

 The Weil pairing (the cup product in  cohomology) $J_2\times J_2\to
 \mu_2$, where $\mu_2=\{+1,-1\}\subseteq \Bbb{C}^*$ induces an isomorphism of groups $W: J_2\to \hat{A}$. The following
 proposition follows from results of ~\cite{bls} (see Section
 ~\ref{california}).
\begin{proposition}\label{july3}
For $\alpha\in J_2=A$,
$$\Pdescend_{\kappa\tensor
\alpha}=P_{\kappa}\tensor L_{W(\alpha)}\in \operatorname\Pic(Y).$$

\end{proposition}
{{\bf Notation:}} Fix an even theta characteristic $\kappa$ for the
rest of this paper. Let $P=\pii^*(P_{\kappa})\in
\Pic(\widetilde{Y})$, $\PP=p^*P_{\kappa}\in
\Pic(\mathcal{M}_{\Spin(r)})$. Denote the descent of $\PP$ to
$M^{\reg}_{\Spin(r)}$ again by $P$. Note that the two definitions of
$P$ are canonically identified under the inclusion
$\widetilde{Y}\subseteq M^{\reg}_{\Spin(r)}$ (using descent theory).
Also note that $p^*P_{\kappa'}$ is isomorphic to $P$, for any theta
characteristic $\kappa'$, the isomorphism is unique upto scalars.

We have a decomposition as $A$-modules:
$$H^0(\widetilde{Y},P)=\oplus_{\chi\in\hat{A}} H^0(Y,\Pdescend_{\kappa}\tensor L_{\chi})$$
\begin{proposition}\label{evensteven} For even $r$,
\begin{enumerate}
\item $H^0(\widetilde{Y},P)$ is $2^{2g}$ dimensional.
\item Each $H^0(Y,\Pdescend_{\kappa}\tensor L_{\chi})$ is one dimensional and spanned by the pfaffian
section of $P_{\kappa\tensor W^{-1}(\chi)}$ corresponding to the
isomorphism in Proposition ~\ref{july3}.
\item The elements $s_{\kappa'}$ in $H^0(\widetilde{Y},P)$ for $\kappa'\in \theta(X)$ form a
basis.
\item The element $s_{\kappa'}$ for  $\kappa'\in \theta(X)$ spans the $\chi=W(\kappa'\tensor\kappa^{-1})$
isotypical subspace of $H^0(\widetilde{Y},P)$.
\end{enumerate}
\end{proposition}
\begin{proof}
$M_{\Spin(r)}-\widetilde{Y}$ has codimension $\geq 2$. Using results
in ~\cite{bls},
\begin{equation}\label{july5}
H^0(\widetilde{Y},P)=H^0(M_{\Spin(r)}^{\reg},P)=H^0(\mathcal{M}_{\Spin(r)},\mathcal{P})=
2^{2g}.
\end{equation}
(for the last equality see ~\cite{ox}.)

 Clearly, the vector space in (2) has at least the (non-zero) pfaffian
section. Since the number of theta-characteristics is $2^{2g}$, (2)
follows from (1). Finally, (3) and (4) are restatements of (2).
\end{proof}

For $r$ odd, we have the following result, whose proof is similar to
that of Proposition ~\ref{evensteven} (recall that our fixed theta
characteristic $\kappa$ is assumed to be even).
\begin{proposition}\label{evensteven} For odd $r$,
\begin{enumerate}
\item $H^0(\widetilde{Y},P)$ is $2^{g-1}(2^g+1)$ dimensional.
\item The elements $s_{\kappa'}$ in $H^0(\widetilde{Y},P)$ for even $\kappa'\in \theta(X)$ form a
basis.
\item The element $s_{\kappa'}$ for even  $\kappa'\in \theta(X)$ spans the $\chi=W(\kappa'\tensor\kappa^{-1})$
isotypical subspace of $H^0(\widetilde{Y},P)$.
\end{enumerate}
\end{proposition}

\section{Hitchin's connection and the geometric Segal-Sugawara tensor}\label{sugawara}
Let $\widetilde{G}$ be a simple, simply connected group and
$M=M_{\widetilde{G}}^{\reg}(X)$ the smooth open subvariety of
$M_{\widetilde{G}}(X)$ parameterizing regularly stable bundles $E$.
Let $\mathcal{M}_g$ denote the moduli stack of smooth and connected
projective algebraic curves of genus $g$, and $X\in \mathcal{M}_g$
as before. The cup product
$$H^1(X,T_X)\tensor H^0(X, \operatorname{ad}(E)\tensor K_X)\to H^1(X,
\operatorname{ad}(E))$$ and the identification $T_E
M=H^1(X,\operatorname{ad}(E))$ defines a (``geometric
Segal-Sugawara'') morphism
$$S:T_X\mathcal{M}_g\to H^0(M, S^2 TM)$$

The group $B$ of principal  $Z(\widetilde{G})$-bundles acts on $M$
and the functoriality of the cup product implies that the  morphism
$S$ has its image in the  subspace of invariants $H^0(M, S^2
TM)^{B}$. {\em We will assume by passing to \'etale covers, that in
any family of curves, the group scheme of principal $B$-bundles
(which sits inside the torsion in a product of Jacobians) has been
trivialized.}

Let $\ml$ be the generating line bundle of the Picard stack of
$\mathcal{M}_{\widetilde{G}}(X)$, which descends to $M$. The action
of $B$ on $M$ may not lift to the pair $(M,\ml)$. For $b\in B$,
$b^*\ml$ is isomorphic to $\ml$ and hence we can form a
 Mumford-theta group $\mg(X)$, a central extension of $B$ by
$\Bbb{C}^*$, which does act on the pair $(M,\ml)$.

In the case $\widetilde{G}=\operatorname{SL}(n)$, it is possible to
identify this Mumford-theta group (the author learned this from M.
Popa, and appears in ~\cite{belkale}). In the case of the odd spin
groups $\Spin(r)$, $r$ odd, the group extension $\mg(X)$ splits,
because replacing $M^{\reg}_{\Spin_r}(X_s)$ by $\widetilde{Y}_s$
(which does not change Picard groups, and isomorphisms of line
bundles over $\widetilde{Y}_s$ extend to $M^{\reg}_{\Spin_r}(X_s)$)
the pfaffian line bundle is pulled back from the (regularly stable)
moduli of odd orthogonal bundles of rank $r$. Here we note that the
center of the odd spin groups is $\Bbb{Z}/2$.

We do not know in general how to identify $\mg(X)$. However suppose
we are given a lifting of the $B$ action on $M$ to an action of a
subgroup $A\subseteq B$ on $(M,\ml)$, where $A$ is the group of
principal $A'$-bundles for some subgroup $A'\subset
Z(\widetilde{G})$. In this setting, it is easy to see by an obvious
generalization of Corollary 5.2 and Lemma 4.1 in ~\cite{belkale},
that
\begin{lemma}\label{rain}
The action of $A$ on $H^0(M,\ml^k)$ preserves the Hitchin connection
as $X$ varies in a family.
\end{lemma}
\begin{remark}
Hitchin's connection is given by ``projective heat operators''. By
averaging over $A$ one can find heat operators invariant under the
action of $A$ (as in ~\cite{r2g2}). Lemma ~\ref{rain} follows
immediately.
\end{remark}

We will now carry out the proof of Theorem ~\ref{phone}. The proof
of Theorem ~\ref{phone2} is similar and hence omitted.

Let us by passing to an open cover in the \'etale topology, assume
that the sheaf of theta characteristics and also the two-torsion in
the Jacobian of the fibers are trivial over $S$. We can form
relative versions of the spaces $\widetilde{Y}$, $Y$ from the
previous discussion. There is an action of $A$ on
$(M_{\Spin(r)}^{\reg},\Pdescend)$, which restricts to the action on
$(\widetilde{Y}_s,{P})$ (because of the codimension estimates).

 Clearly, by the fiberwise
equality ~\eqref{july5},
$$H^0(\mathcal{M}_{\Spin(r)}({X}_s),
\mathcal{P})=H^0(\widetilde{Y}_s,{P}).$$ We have an action of the
(trivial group scheme) $A=J_2$, corresponding to
$A'=\ker(\Spin(r)\to \SO(r))$ on the right hand side. This action
preserves the Hitchin connection (by ~\eqref{july5} and Lemma
~\ref{rain}). Given this it is easy to finish the argument. The
isotypical components of the action of $A$ are preserved by the
Hitchin connection (this is obvious if we choose a $A$-invariant
heat operator). In particular each of the isotypical spaces, each
spanned by some $s_{\kappa'}$ is preserved by the Hitchin
connection.

\section{The proof of Proposition ~\ref{july3}}\label{california}

The arguments in this section are taken from ~\cite{bls} and Section
5.3 of ~\cite{L}. Fix a point $\point\in X$ and a formal coordinate
$z$ at $\point$. For ease of notation let $\widetilde{G}=\Spin(r)$
and $G=\SO(r)$, $LG=G(\Bbb{C}((z)))$, $L_XG=
G(\mathcal{O}(X-\point))$, $L^+G=G(\Bbb{C}[[z]])$ (similarly define
$L\widetilde{G}$, $L_X G$ and $L^+\widetilde{G}$). We have two
``infinite'' Grassmannians
$$\mq_{G}= LG/L^+G,\  \mq_{\widetilde{G}}= L\widetilde{G}/L^+\widetilde{G}$$

The space $\mq_G$ (similarly $\mq_{\widetilde{G}}$) parameterizes
isomorphism classes of principal $G$-bundles equipped with a
trivialization on $X-\{\point\}$.

 It is known from ~\cite{bls}, that the neutral component $\mq_G^o$
of $\mq_{G}$ is canonically isomorphic to $\mq_{\widetilde{G}}$.
Hence a $G$-bundle (in the neutral component of $\mathcal{M}_G$)
trivialized on the complement of $\point$ has a canonical
$\widetilde{G}$-structure. It is also known that $L_X G$ is
contained in the neutral component of $LG$. Finally, one has the
stack-theoretic uniformization theorems \cite{v2,ls,bls}
$$\mathcal{M}_G = L_X G\backslash \mq_G, \mathcal{M}_{\widetilde{G}} = L_X\widetilde{G}\backslash
\mq_{\widetilde{G}}.$$

Let us show that $L_X G$ acts on $\widetilde{Y}$. Let $P\in
\widetilde{Y}$ and $\beta\in L_XG$.  Represent $P$ as the image of a
point $q\in \mq_{\widetilde{G}}$ and hence as a point of $\mq_G$.
Clearly $L_X G$ acts on $\mq_G$ preserving the connected components.
Therefore $\beta q$ gives a new point of $\mq_{\widetilde{G}}$, and
hence a new point of $\widetilde{Y}$. In fact this action of $L_X G$
factors through the quotient by image of $L_X \widetilde{G}$. The
quotient $L_X G/i(L_X\widetilde{G})$ is naturally isomorphic to
$J_2=A$ (see Lemma 1.2 in ~\cite{bls}), and this action of $A$ on
$\widetilde{Y}$ coincides with the natural Galois action of $A$ on
$\widetilde{Y}$ (see Section ~\ref{anti}). There is another way to
describe this action. There is a natural map $L_X G\to
L\widetilde{G}/\pi_1(G)$ (both sides sit naturally in $LG$). Through
this map $L_X G$ acts on $\widetilde{Y}$, and it is easy to see that
it coincides with the action above (the natural map
$\mq_{\widetilde{G}}\to \mq_G$ is equivariant for the map of groups
$L\widetilde{G}\to LG$).

 In Section 5 of ~\cite{bls}, an injective
homomorphism $\lambda: \hat{A}\to \Pic(\mathcal{M}_{G})$ is
constructed and it is shown that as line bundles on
$\mathcal{M}_{G}$, $\PP_{\kappa\tensor \alpha}\tensor
\PP_{\kappa}^{-1}$ equals $\lambda(W(\alpha))$ (see the proof of
Proposition 5.2 in ~\cite{bls}). We claim that the descent of
$\lambda(\chi)$ to $Y$ equals $L_{\chi}$ for $\chi\in\hat{A}$. This
would prove Proposition ~\ref{july3}.

For simplicity, we will work in the classical topology over $Y$
(which is sufficient for our purposes, because of the codimension
conditions). In fact, it is easy to replace the argument by an
analogous argument in the \'etale topology, and prove Proposition
~\ref{july3} in the algebraic category. Let us first recall our
construction of $L_{\chi}$. Cover $Y$ by (analytic) open subsets
$U_i$ and choose a lifting $U_i\to \widetilde{Y}$. On overlaps
$U_i\cap U_j$, the two maps  differ by a section of $A$. Hence a
character $\chi$ of $A$ gives the patching functions for a line
bundle on $Y$ (which coincides with $L_{\chi}$, note that
$\chi=\chi^{-1}$).

We will now realize this construction by making loop group choices.
This is then easily seen to be the construction in ~\cite{bls}:
Refine the cover $U_i$ and on each $U_i$ choose a local universal
bundle $Q_i$ (this is possible using $Y\subseteq
M^{\reg}_{\SO_r}(X)$) and a trivialization of $Q_i$ on the
complement of $\point$. This gives $Q_i$ a
$\widetilde{G}$-structure, and hence we obtain liftings $U_i\to
\widetilde{Y}$. On overlaps $U_i\cap U_j$, the different
trivializations give a class in $L_X G/Z(G)$. Therefore any
character $\chi$ of $L_X G/Z(G) $, produces a line bundle  on $Y$.
Any such character is necessarily trivial on the image of
$L_X\widetilde{G}$, and factors through the quotient $L_X
G/i(L_X\widetilde{G})=A$ where $i:L_X\widetilde{G}\to L_X G$ (note
that the center of $\widetilde{G}$ surjects on to the center of
$G$). By the basic compatibility verification in Section
~\ref{anti}, the proof of our claim is complete.

\section{Application to the Symplectic strange
duality}\label{antimatter}
 Let us first recall the set up    of the
symplectic strange duality from ~\cite{bea}. Let $M_{\Sp(2n)}$
denote the moduli space of vector bundles on $X$ of rank $2n$,
equipped with a non-degenerate symplectic form (with values in
$\mathcal{O}_X$). In fact, $M_{\Sp(2n)}$ is the moduli space of
principal $\Sp(2n)$-bundles on $X$. Let  $\ml$ be the positive
generator of the Picard group of $M_{\Sp(2n)}$. We can take $\ml$ to
the determinant of cohomology of the tautological bundle tensored
with a degree $g-1$ line bundle on $X$ (this makes good sense on the
moduli stack, and descends to the moduli space).

Similarly let $M'_{\Sp(2m)}$ denote the moduli space of vector
bundles on $X$ of rank $2m$, equipped with a non-degenerate
symplectic form with values in $K_X$ (therefore the underlying
degree of the vector bundles is $2m(g-1)$). A choice of a theta
characteristic $\kappa$ gives an isomorphism $M_{\Sp(2m)}\to
M'_{\Sp(2m)}$. Let $\ml$ again denote the positive generator of the
Picard group of $M'_{\Sp(2m)}$. Note that for both $M_{\Sp(2n)}$ and
$M'_{\Sp(2m)}$, the global sections of powers of $\ml$ over the
corresponding moduli stack, coincides with the global sections over
the moduli space.

On $M_{\Sp(2n)}\times M'_{\Sp(2m)}$, there is a natural Cartier
divisor $\Delta$ of the line bundle $\ml^m\boxtimes\ml^n$, such that
$2\Delta$ is the theta section of the determinant of cohomology of
the tensor product. The non-zeroness of this divisor has been shown
by Beauville ~\cite{bea}. Therefore one finds a non-zero
homomorphism, conjectured by Beauville to be an isomorphism
\begin{equation}\label{SD}
H^0({M}'_{\Sp(2m)}(X),\mathcal{L}^n)^*\to
H^0({M}_{\Sp(2n)}(X),\mathcal{L}^m)
\end{equation}
Said in a different way, the divisor on the product of the
moduli-stacks $\mathcal{M}_{\Sp(2n)}\times \mathcal{M}_{\Sp(2m)}$ is
the  pull back of the pfaffian section $s_{\kappa}$ on
$\mathcal{M}_{\Spin(4mn)}$ of the line bundle $\PP$. That is, the
image of $s_{\kappa}$ under the map
\begin{equation}\label{beauville}
H^0(\mathcal{M}_{\Spin(4mn)}(X),\mathcal{P})\to
H^0(\mathcal{M}_{\Sp(2m)}(X),\mathcal{L}^n)\times
H^0(\mathcal{M}_{\Sp(2n)}(X),\mathcal{L}^m)
\end{equation}
It is known that the map ~\eqref{beauville} is projectively flat
(see ~\cite{nt} and ~\cite{belkale}). Therefore, by Theorem
~\ref{phone}, we see that the map ~\eqref{SD} is a projectively flat
map after making the identification $M_{\Sp(2m)}\to M'_{\Sp(2m)}$.
We therefore obtain the following corollary to Theorem ~\ref{phone}:
\begin{corollary}\label{dual} The homomorphism ~\eqref{SD} has constant rank as $X$ varies over the moduli space of curves $M_g$.

 \end{corollary}
\subsection{Relations to the strange duality for vector bundles}\label{etrange} Consider the case $n=1$ and (for technical
reasons) $g>2$.  By the above results, the local system on the
moduli of curves with a choice of theta characteristic, given by
$H^0({M}_{\Sp(2m)}(X),\mathcal{L})$ is naturally (projectively) dual
to the local  system with fibers
$$H^0({M}_{\Sp(2)},\mathcal{L}^m)=H^0({M}_{\SL(2)},\mathcal{L}^m).$$
Using the  ${\SL}(2)$-${\operatorname{GL}}(m)$ strange duality, and
its flatness
 ~\cite{L2,a,b1,mo,belkale} we find that the latter space is naturally
dual, preserving connections to
$H^0(M_{\operatorname{GL}(m)}(0),\mathcal{L}^2)$, where
$M_{\operatorname{GL}(m)}(0)$ is the moduli space of semi-stable
degree $0$ and rank $m$ vector bundles\footnote{Note that we have
made a choice of a theta characteristic on $X$, and the line bundle
on $\ml$ on $M_{\operatorname{GL}(m)}(0)$ is the determinant of
cohomology of ``the tautological bundle'' $\tensor \kappa$, which
descends from the corresponding stack. In fact $\ml^2$ does not
depend upon $\kappa$ but its strange duality with
$H^0({M}_{\SL(2)},\mathcal{L}^m)$ does depend on $\kappa$.} on $X$.
Actually, there is a natural embedding
$M_{\operatorname{GL}(m)}(0)\subseteq {M}_{\Sp(2m)}$ which pulls
back $\mathcal{L}$ to $\mathcal{L}^2$, and is consistent with the
above identifications. Therefore,  the natural
 map $H^0({M}_{\Sp(2m)},\mathcal{L})\to
 H^0(M_{\operatorname{GL}(m)}(0),\mathcal{L}^2)$
is an isomorphism, preserving connections.

 Note that $\operatorname{GL}(m)\subseteq \Sp(2m)$
appears as a conformal embedding in the tables of conformal
embeddings, but the author does not know how to use this to directly
prove  that the natural map from $H^0({M}_{\Sp(2m)},\mathcal{L})$ to
$H^0(M_{\operatorname{GL}(m)}(0),\mathcal{L}^2)$ preserves
connections (the problem is the non-semisimplicity of
$\operatorname{GL}$).

\section{A verification of compatibility}\label{anti}

Let $G$ be a semisimple algebraic group, with universal cover
$\widetilde{G}$. We have a basic central extension
$$1\to \pi_1(G)\to \widetilde{G}\leto{\pi} G\to 1$$

Let $X$ be a smooth projective curve as before and $\point$ a point
on it. Set $X^*=X-\{\point\}$, and consider a map $\phi:X^*\to G$,
or an element $\phi\in L_X G$ using our earlier notation. We find by
base change a cover $\widetilde{X}^*$ of $X^*$ which fits into a
cartesian diagram
\begin{equation}\label{digram}
\xymatrix{ \widetilde{X}^*\ar[r]^{\tilde{\phi}}\ar[d]^{\pi'} &  \widetilde{G}\ar[d]^{\pi}\\
{X}^*\ar[r]^{{\phi}} &{G}}
\end{equation}

Now unramified abelian covers of $X^*$ extend to unramified abelian
covers of $X$. Therefore we can extend $\pi'$ to a cover
$\pi':\widetilde{X}\to X$ and thus obtain  a principal
$\pi_1(G)$-bundle $\alpha$ on $X$ in the \'{e}tale topology. Given a
principal $\widetilde{G}$-bundle $Q$ on $X$ we can
  obtain a new bundle $Q_1$ on $X$ whose sheaf of sections is for an
  open subset $U$ of $X$, sections of the pull back of $Q$ over the
  inverse image of $U$ which twist by the image of $\pi_1(G)$ in
  $\widetilde{G}$ upon the action of the covering group $\pi_1(G)$. It
  is easy to see that $Q_1$ is the same as $Q\times_{\pi_1(G)}
  \alpha$ (this leads to the Galois action of $\alpha$ on
  $M_{\widetilde{G}}$).

On the other hand, given $Q$ we have another construction of a
principal $\widetilde{G}$-bundle on $X$. There is a natural map $L_X
G\to L\widetilde{G}/\pi_1(G)$. To do this pick a point
$\tilde{\point}$ over $\point$ and a coordinate $z$ on $X$ at
$\point$. Since $\pi'$ is \'etale, $z$ lifts to a coordinate  on
$\widetilde{X}$ near $\tilde{\point}$. The map $\tilde{\phi}$
therefore gives us an element $\psi\in L\widetilde{G}$, which
normalizes $L_X{\widetilde{G}}$ (by descent theory) and hence left
multiplication by $\psi$ gives a principal $\widetilde{G}$-bundle
$Q_2$ on $X$. We contend that $Q_1$ and $Q_2$ are isomorphic.

Let $s$ be a section of $Q$ over $X^*$, clearly $\tilde{\phi}s$
gives a section of $Q_1$ over $X^*$. Also a section of $Q_1$ over
$D$ a formal neighborhood of $\point$ and the choice of
$\tilde{\point}$ over $\point$, determines a section of $Q_1$ over
$D$. The new patching function for $Q_1$ (in the punctured disc
around ${\point}$) is given by the image $\psi$ of $\tilde{\phi}$
times the patching function of $Q$, hence $Q_1$ is isomorphic to
$Q_2$.

The following diagram (easily seen to be commutative) is useful in
studying the various maps, where the vertical arrow is the map
$\phi\to\tilde{\phi}$ as above:
\begin{equation}
\xymatrix{ L\widetilde{G}/\pi_1(G)\ar[r] &  LG\\
L_X G\ar[u]\ar[ur]}
\end{equation}

\subsection{Action of the center}\label{consequence} The above discussion has the following
interesting consequence (take $G=\widetilde{G}/Z(\widetilde{G})$):
The action of a principal $Z(\tilde{G})-$bundle on the set of
isomorphism classes of principal $\widetilde{G}$-bundles on $X$,
lifts to left multiplication by an element of $L\widetilde{G}$ on  $
\mq_{\widetilde{G}}$.

\bibliographystyle{plain}
\def\noopsort#1{}

\end{document}